\documentclass[12pt]{article}
\usepackage{amssymb,amsfonts}
\usepackage{graphicx}
\usepackage{amscd}
\newtheorem{theorem}{Theorem}
\newtheorem{lemma}{Lemma}
\newtheorem{prop}{Proposition}
\newtheorem{corollary}{Corollary}
\newtheorem{definition}{Definition}
\title{Well-posed Infinite Horizon Variational Problems on a Compact Manifold}
\author{A. Agrachev \thanks{SISSA, Trieste and MIAN, Moscow}}
\begin{document}
\maketitle

\begin{abstract}
We give an effective sufficient condition for a variational
problem with infinite horizon on a compact Riemannian manifold $M$
to admit a smooth optimal synthesis, i.\,e. a smooth dynamical
system on $M$ whose positive semi-trajectories are solutions to
the problem. To realize the synthesis we construct a
well-projected to $M$ invariant Lagrange submanifold of the
extremals' flow in the cotangent bundle $T^*M$. The construction
uses the curvature of the flow in the cotangent bundle and some
ideas of hyperbolic dynamics.
\end{abstract}
\section{Introduction}
This paper is dedicated to the 100th anniversary of Lev Semenovich
Pontryagin.

Let $M$ be a smooth compact $n$-dimensional Riemannian manifold.
Given \linebreak $\alpha\ge 0$, we denote by $\Omega^\alpha$ the set
of all absolutely continuous curves \linebreak
$\gamma:[0,+\infty)\to M$ such that the integral
$\int\limits_0^{+\infty}e^{-\alpha t}|\dot\gamma(t)|^2\,dt$
converges. Here $|\dot\gamma(t)|$ is the Riemannian length of the
tangent vector $\dot\gamma(t)\in T_{\gamma(t)}M$. Given $q\in M$, we
set $\Omega^\alpha_q=\{\gamma\in\Omega^\alpha : \gamma(0)=q\}$.

Let $U:M\to\mathbb R$ be a smooth function. We set
$$
\mathfrak I_\alpha(\gamma)= \int\limits_0^{+\infty}e^{-\alpha
t}\left(\frac 12|\dot\gamma(t)|^2-U(\gamma(t)\right)\,dt
$$
and try to minimize $\mathfrak I_\alpha$ on $\Omega^\alpha_q$. Given
$q_0\in M$, we say that $\gamma_0\in\Omega^\alpha_{q_0}$ is a {\it
minimizer} if $ \mathfrak I_\alpha(\gamma_0)=\min\{\mathfrak
I_\alpha(\gamma):\gamma\in\Omega^\alpha_{q_0}\}. $ We say that
$\mathfrak I_\alpha$ defines a well-posed variational problem or
{\it admits smooth optimal synthesis} if for any $q\in M$ there
exists a unique minimizer $\gamma_q$ and the map $(q,t)\mapsto
\dot\gamma_q(t),\ q\in M,\,t\ge 0$, is of class $C^1$.

The functional $\mathfrak I_0$ is simply the action of a mechanical
system on the Riemannian manifold $M$ with potential energy $U$. If
$\alpha>0$, then $\mathfrak I_\alpha$ is the {\it discounted action}
with the {\it discount factor} $\alpha$.

It is not hard to show that $\mathfrak I_0$ does not define a
well-posed variational problem for generic $U$. In this paper we
prove that $\mathfrak I_\alpha$ defines a well-posed problem for
all sufficiently big $\alpha$. Moreover, we give the effective
sharp estimate for critical $\alpha$.

The functional $\mathfrak I_\alpha$ admits smooth optimal
synthesis if and only if there exists a complete vector field $V$
of class $C^1$ on the manifold $M$ such that any positive
semi-trajectory $\gamma$ of the dynamical system $\dot q= V(q)$ is
a unique minimizer of $\mathfrak I_\alpha$ on
$\Omega^\alpha_{\gamma(0)}$. Indeed, assume that $\mathfrak
I_\alpha(\gamma(\cdot))=\min\{\mathfrak
I_\alpha(q(\cdot)):q(\cdot)\in\Omega^\alpha_{q_0}\}$, then
$\mathfrak I_\alpha(\gamma(s+\cdot))= \min\{\mathfrak
I_\alpha(q(\cdot)):q(\cdot)\in\Omega^\alpha_{\gamma(s)}\}$ for any
$s\ge 0$. Define $V(q_0)=\dot\gamma_{q_0}(0)$; then $q_0\mapsto
V(q_0)$ is a vector field of class $C^1$ on $M$ and
$\dot\gamma_{q_0}(s)=V(\gamma(s)),\ \forall s\ge 0$.

How to characterize the minimizers? If $\gamma$ is a minimizer
then, obviously, $\gamma$ minimizes the functional
$\int\limits_0^Te^{-\alpha t}\left(\frac 12|\dot
q(t)|^2-U(q(t)\right)\,dt$ on the space of all $q(\cdot)\in
H^1([0,T],M)$ such that $q(0)=q_0,\,q(T)=\gamma(T)$ for any $T>0$.
Hence any solution of our infinite horizon variational problem
satisfies the Euler--Lagrange equation of the classical finite
horizon variational problem.

Let us consider a simple example. Let  $M=S^1=\mathbb
R/2\pi\mathbb Z$ and $U(\theta)=\cos\theta,\ \theta\in S^1$, the
potential energy of the mathematical pendulum. The Euler--Lagrange
equation has the form: $\ddot\theta=\sin\theta+\alpha\dot\theta$.
Write it as a system:
$$
\left\{\begin{array}{rcl} \dot\theta&=&\xi\\
\dot \xi&=&\sin\theta+\alpha\xi\\ \end{array}\right.. \eqno (1)
$$
This system has 2 equilibrium points: $(\theta,\xi)=(0,0)$ and
$(\theta,\xi)=(\pi,0)$. The equilibrium $(0,0)$ is a saddle for any
$\alpha\ge 0$. The equilibrium $(\pi,0)$ is a center for $\alpha=0$
(see Fig.~1), unstable focus for $0<\frac{\alpha^2}4<1$ (Fig.~2),
and unstable node for $\frac{\alpha^2}4>1$ (Fig.~3).

\begin{figure}
\hspace{3cm} \begin{minipage}[b]{6.5cm}
   \centering
   \includegraphics[width=7cm]{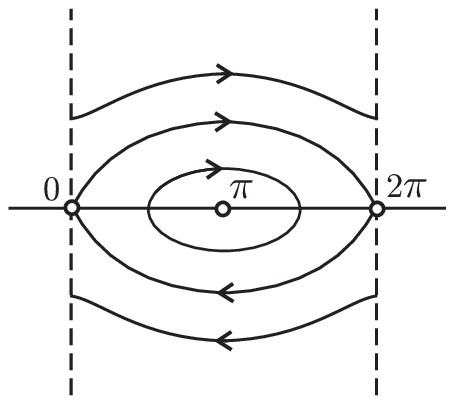}
   \caption{}
 \end{minipage}
 \end{figure}

\begin{figure}
 \begin{minipage}[c]{6.5cm}
  \centering
   \includegraphics[width=6cm]{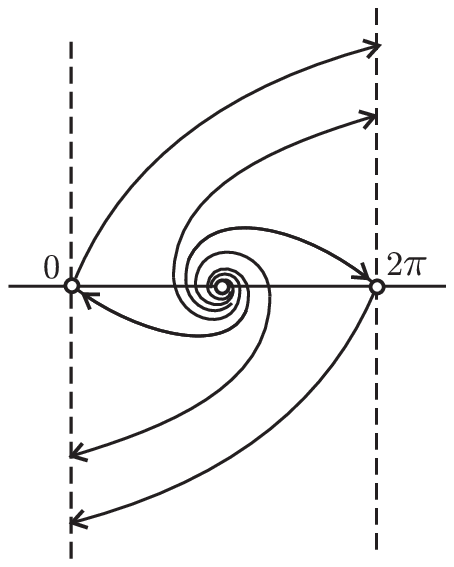}
   \caption{}
 \end{minipage}
 \ \hspace{2mm} \hspace{3mm} \
 \begin{minipage}[c]{6.5cm}
  \centering
  \includegraphics[width=6cm]{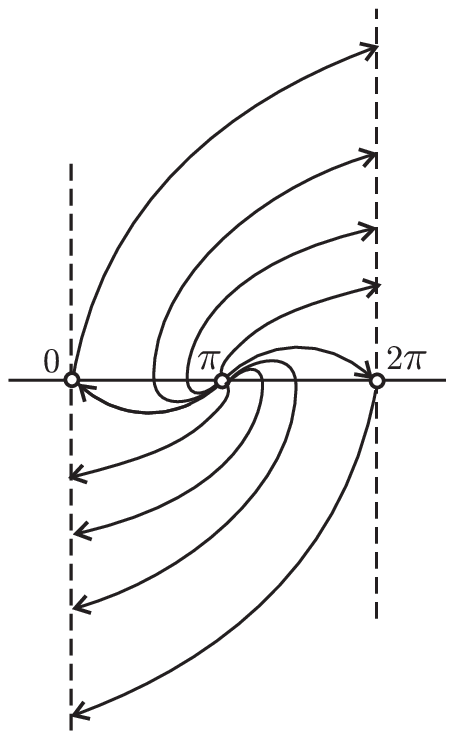}
   \caption{}
 \end{minipage}
\end{figure}

A solution $\theta(\cdot)$ of the Euler--Lagrange equation belongs
to $\Omega^\alpha$ if and only if $(\xi(\cdot),\theta(\cdot))$ is
an equilibrium or a part of the stable submanifold of the saddle.
The saddle and node equilibria are minimizers while the focus and
center are not. If $\frac{\alpha^2}4<1$ then $\mathfrak I_\alpha$
does not admit smooth optimal synthesis; if $\frac{\alpha^2}4>1$
then it does admit smooth optimal synthesis. The corresponding to
the minimizers trajectories of system (1) in both cases are shown
on Figure~4.

\begin{figure}
\hspace{3cm} \begin{minipage}[b]{6.5cm}
   \centering
   \includegraphics[width=8cm]{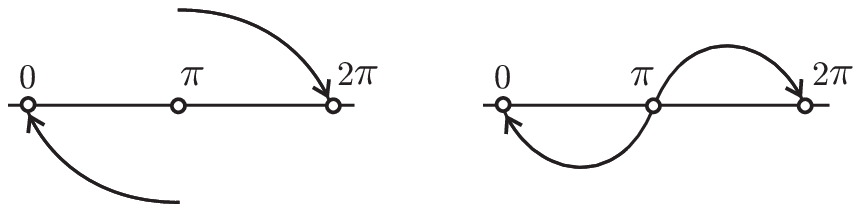}
   \caption{}
 \end{minipage}
 \end{figure}

In order to formulate our main result we need the following
definition.
\begin{definition} The curvature of the Hamiltonian $H$ at $z\in
T^*_qM$ is a self-adjoint linear operator $R_z^H:T^*_qM\to T^*_qM$
defined by the formula
$$
R^H_z\zeta=\mathfrak R(\zeta,z)z+(\nabla^2_qU)\zeta,\quad \zeta\in
T^*_qM,
$$
where $\nabla$ is the covariant derivativation and $\mathfrak R$
the Riemannian curvature.
\end{definition}

Given a self-adjoint linear operator $A$ on a Euclidean space and
a constant $a\in\mathbb R$, the relation $A<aI$ means that all
eigenvalues of $A$ are smaller than $a$.
%Let $\varpi:\tilde M\to M$ be a covering of $M$. The pullback
%$\tilde J_\alpha$ of the functional $J_\alpha$ is defined on the
%curves in $M$ by the formula $\tilde
%J_\alpha(\gamma)=J_\alpha(\varpi\circ\gamma)$.
Now we state the main result of this paper:

\begin{theorem}
Let $R^H_z<\frac{\alpha^2}4I$ for any $z\in T^*M$ such that
$H(z)\le\max\limits_{q\in M}U(q)$; then
%there exists a finite covering $\varpi:\tilde M\to M$ such that
$\mathfrak I_\alpha$ admits smooth optimal synthesis.
\end{theorem}

\begin{corollary} Assume that
sectional curvature of $M$ does not exceed $r\ge 0$ and
$$
\nabla^2_qU<\left(\frac{\alpha^2}4-2r(\max U-\min U)\right)I,\quad
\forall q\in M.
$$
%implies the existence of a finite covering $\varpi:\tilde M\to M$
%such that
Then $\mathfrak I_\alpha$ admits smooth optimal synthesis.
\end{corollary}

In the next section we use symplectic language to characterize
extremals of the variational problem and formulate a more detailed
version of the main result, including its Hamilton--Jacobi
interpretation. The discount factor appears as a friction
coefficient in the extremals' equation and serves as a ``smoothing
factor" in the Hamilton--Jacobi setting.

In Section 3 we prove the ``partial hyperbolicity" of the set
filled by the extremals; this is a principal step in the proof of
the main result. A good source for the theory of partially
hyperbolic systems is book \cite{Pe}. It seems that this concept
is really relevant to our subject: sufficiently strong isotropic
friction almost automatically leads to the partial hyperbolicity
and this fact concerns much more general systems than ones studied
in this paper.

In Section 4 we complete the proof of the main result: optimal
synthesis is obtained as a limit of solutions of the horizon
$\tau>0$ problem as $\tau\longrightarrow+\infty$.

In this paper we assume that the manifold $M$ is compact. This
assumption can be substituted by requirements on the asymptotic
behavior of the Riemannian metric and potential at infinity. In
particular, the result remains valid for the ``asymptotically
flat" metrics and ``asymptotically quadratic" potentials. This and
other generalizations (more general Lagrangians, problems with
nonholonomic constraints) could be subjects of forthcoming papers.

\medskip\noindent{\sl Acknowledgment.} I am grateful to
Dmitriy Treschev for the stimulating discussion.

\section{Symplectic setting}

Let $\sigma$ be the standard symplectic form on $T^*M$,
$\sigma=ds$, where $s$ is the Liouville form:
$s_z=z\circ\pi_*\bigr|_{T_z(T^*M)},\ \forall z\in T^*M$, and
$\pi:T^*M\to M$ is the standard projection.

 The Hamiltonian (or the energy function) $H:T^*M\to\mathbb R$
associated to the Lagrangian $\frac
12|\dot\gamma(t)|^2-U(\gamma(t))$ is defined by the formula
$H(z)=\frac 12|z|^2+U(\pi(z)),\ \forall z\in T^*M$, where
$|z|=\max\{\langle z,v\rangle:v\in T_{\pi(z)}M,\ |v|=1\}$ is the
dual norm to the norm defined by the Riemannian structure.
Actually, the Riemannian structure is a self-adjoint isomorphism
of $TM$ and $T^*M$ and in this paper we freely use without special
mentioning the provided by this isomorphism identification of
tangent and cotangent vectors.

The Legendre transformation of the time-varying Lagrangian
\linebreak $e^{-\alpha t}\left(\frac 12|\dot q|^2-U(q)\right)$ has
the form $e^{-\alpha t}H(e^{\alpha t}z)$, hence solutions of the
Euler--Lagrange equations are exactly projections to $M$ of the
trajectories of the Hamiltonian system on $T^*M$
$$
\dot z=e^{-\alpha t}\vec H(e^{\alpha t}z), \eqno (2)
$$
where $\vec H$ is the Hamiltonian vector field associated to $H$;
the field $\vec H$ is defined by the identity
$dH=\sigma(\cdot,\vec H)$. Recall that $z$ is a point of the
$2n$-dimensional manifold $T^*M$; local coordinates on $M$ induce
local trivialization of $T^*M$ so that $z$ splits in two
$n$-dimensional vectors, $z=(p,q)$, and $\pi:(p,q)\mapsto q$. Then
$s=\sum\limits_{i=1}^np^idq^i,\
\sigma=\sum\limits_{i=1}^ndp^i\wedge dq^i$ and system (2) takes
the form:
$$
\left\{\begin{array}{rcl} \dot q&=&\frac{\partial H}{\partial p}(e^{\alpha t}p,q)\\
\dot p&=&-e^{-\alpha t}\frac{\partial H}{\partial q}(e^{\alpha
t}p,q)\\ \end{array}\right..
$$
The time-varying change of variables $\xi=e^{\alpha t}p$ gives the
system
$$
\left\{\begin{array}{rcl} \dot q&=&\frac{\partial H}{\partial\xi }(\xi,q)\\
\dot \xi&=&-\frac{\partial H}{\partial q}(\xi,q)+\alpha\xi\\
\end{array}\right.
$$
or in the coordinate free form:
$$
\dot\zeta=\vec H(\zeta)+\alpha e(\zeta), \eqno (3)
$$
where $\zeta(t)=e^{\alpha t}z(t)$ and $e$ is the vertical {\it Euler
vector field} of the vector bundle $\pi:T^*M\to M$.

We set $\mathfrak h_\alpha=\vec H+\alpha e$ and denote by
$e^{t\mathfrak h_\alpha},\ t\in\mathbb R$, the flow on $T^*M$
generated by the vector field $\mathfrak h_\alpha$. It is easy to
see that
$$
\left(e^{t\mathfrak h_\alpha}\right)^*s=e^{\alpha
t}s+da_t,\ \forall t\in\mathbb R,
$$
where $a_t$ is a smooth function on $T^*M,\ t\in\mathbb R$. Indeed, let
$s_t=\left(e^{t\mathfrak h_\alpha}\right)^*s$; then
$$
\frac d{dt}s_t=\left(e^{t\mathfrak h_\alpha}\right)^*L_{\mathfrak h_\alpha}s=
\left(e^{t\mathfrak h_\alpha}\right)^*d(i_{\vec H}s-H)+\alpha s_t=
db_t+\alpha s_t,
$$
where $b_t=(i_{\vec H}s-H)\circ e^{t\mathfrak h_\alpha}$. Hence $s_t=e^{\alpha t}s_0+
d\int\limits_0^te^{\alpha(t-\tau)}b_\tau\,d\tau$.

Recall that Lagrange subspaces of $T_z(T^*M)$ are $n$-dimensional
isotropic subspaces of the form $\sigma_z,\ z\in T^*M$, and
Lagrange submanifolds of $T^*M$ are $n$-dimensional submanifolds
whose tangent spaces in all points are Lagrange subspaces.
In other words, an $n$-dimensional submanifold $\mathcal L\subset T^*M$
is a Lagrange submanifold if and only if $s\bigr|_{\mathcal L}$ is a closed form.
We say that $\mathcal L$ is an {\it exact} Lagrange submanifold if $s\bigr|_{\mathcal L}$
is an exact form.
We obtain that $e^{t\mathfrak h_\alpha}$ transforms (exact) Lagrange
submanifolds of $T^*M$ into (exact) Lagrange submanifolds.

More notations. Let $0_q$ be the origin of the space $T^*_qM$. We
set
$$
C_H=\{0_q:d_qU=0,\ q\in M\},
$$
the set of equilibrium points of system (3), and
$$
B_H=\{z\in T^*M : H(z)\le\max U\}.
$$

In what follows we tacitly assume that $\alpha>0$. Note that
$\frac d{dt}H(\zeta(t))=\alpha|\zeta(t)|^2$ for any solution
$\zeta$ of (3). Hence $H$ is strictly monotone increasing along
any solution that is not an equilibrium.

\begin{lemma} Let $\zeta(t)=e^{t\mathfrak h_{\alpha}}(\zeta(0))$
be a trajectory of system (3). The following statements are
equivalent:
\begin{enumerate}
\item $\pi\circ\zeta(\cdot)\in\Omega^\alpha$.
\item $\zeta(\cdot)$ is a bounded curve in $T^*M$.
\item $\zeta(t)\in B_H,\ \forall t\in\mathbb R$.
\end{enumerate}
\end{lemma}
{\bf Proof.}
%The statements {\it 2} and {\it 4} are equivalent due to the
%monotonicity of $H$ along the trajectories.
The implications
{\it3} $\Rightarrow$ {\it 2} $\Rightarrow$ {\it 1} are obvious. Let us prove
that {\it 1} implies {\it 3}.

Assume that the trajectory $\zeta(t),\ t\in\mathbb R$, is not
contained in $B_H$. Then we may assume, without lack of generality,
that $\zeta(0)\notin B_H$. In other words, $H(\zeta(0))-\max U>0$.
Set $\mu(t)=H(\zeta(t))-\max U$, then $\mu(0)>0$. We have:
$$
\dot\mu(t)=\frac
d{dt}H(\zeta(t))=\alpha|\zeta(t)|^2=2\alpha\left(H(\zeta(t))-U(\pi(\zeta(t)))\right)
\ge 2\alpha\mu(t).
$$
Hence $\mu(t)\ge e^{2\alpha t}\mu(0),\ \forall t\ge 0$. Then
$\frac 12|\zeta(t)\|^2\ge e^{2\alpha t}\mu(0)$ and $e^{-\alpha
t}|\zeta(t)|^2\longrightarrow+\infty$ as
$t\longrightarrow+\infty$. Hence $\zeta(\cdot)\notin\Omega^\alpha.
\qquad \square$

\begin{definition}
Extremal locus of the functional $\mathfrak I_\alpha$ is the
subset $\mathcal E_\alpha$ of $T^*M$ filled by those trajectories
of the flow $e^{t\mathfrak h_\alpha}$ that satisfy one of the
conditions 1--3 of Lemma~1 (and hence all these conditions).
\end{definition}

\begin{corollary} $\mathcal E_\alpha$ is a compact invariant
subset of the flow $e^{t\mathfrak h_{\alpha}},\ t\in\mathbb R$.
\end{corollary}

The potential energy $U$ is a {\it Morse function} if $U$ has only
non-degenerate critical points, i.\,e. the Hessian of $U$ at $q\in
M$ is a non-degenerate quadratic form for any $q$ such that
$d_qU=0$.

\begin{lemma}
Assume that $U$ is a Morse function and $\zeta(t)=e^{t\mathfrak
h_{\alpha}}(\zeta(0)),\ t\in\mathbb R$. The curve $\zeta(\cdot)$
is bounded if and only if there exists
$\lim\limits_{t\to+\infty}\zeta(t)\,\in\,C_H.$
\end{lemma}
This lemma is a simple corollary of the monotonicity of $H$ along
trajectories of the flow $e^{t\mathfrak h_{\alpha}}.
\qquad\square$

Next statement is an expanded version of Theorem~1; it explains
the way the smooth optimal synthesis is constructed.

\begin{theorem} If $R_z^H<\frac{\alpha^2}4I,\ \forall z\in B_H$, then
the flow $e^{t\mathfrak h_\alpha}$ has an invariant exact Lagrange
submanifold $\Psi\subset T^*M$ of class $C^1$ such that
$\pi\bigr|_{\Psi}$ is a diffeomorphism of $\Psi$ on $M$ and the
minimizers are exactly positive semi-trajectories of the dynamical
system $\dot q=V(q)$, where $V=\pi_*\left(\mathfrak
h_\alpha\bigr|_{\Psi}\right)$. Moreover,
$C_H\subset\Psi\subset\mathcal E_\alpha$; if $U$ is a Morse function
then $\Psi=\mathcal E_\alpha$.
\end{theorem}

\noindent{\bf Remark 1.}  Submanifold $\Psi\subset T^*M$ in Theorem~2 is the graph
of an exact 1-form on $M$, i.\,e. $\Psi=\{d_qu : q\in  M\}$, where
$u$ is a $C^2$-function on $M$. Then $u$ is a solution of the modified
Hamilton--Jacobi equation:
$$
H(du)-\alpha u=const.
$$
Indeed, let $q\in M$; then $T_{d_qu}\Psi$ is a Lagrange subspace of
$T_{d_qu}(T^*M)$ and $\mathfrak h_\alpha\in T_{d_qu}\Psi$. Hence
$$
0=\sigma(\xi,\mathfrak h_\alpha)=\sigma(\xi,\vec H)+\alpha\sigma(\xi,e)=
\langle d_{d_qu}H-\alpha d_qu,\xi\rangle,\quad \forall\xi\in T_{d_qu}\Psi.
$$
We obtain that $d(H(du)-\alpha u)=0$. The choice of $const$ is at
our disposal. If we set $const=0$, then $-u(q)=\min\left(\mathfrak
I_\alpha\bigr|_{\Omega^\alpha_q}\right),\ \forall q\in M$. Indeed,
$\min\left(\mathfrak
I_\alpha\bigr|_{\Omega^\alpha_q}\right)=\mathfrak
I_\alpha(\gamma)$, where $\gamma(t)=\pi\circ e^{t\mathfrak
h_\alpha}(q)$. We have: $ \mathfrak
I_\alpha(\gamma)=\int\limits_0^\infty e^{-\alpha t}\left(\langle
d_{\gamma(t)}u,\dot\gamma(t)\rangle-H(d_{\gamma(t)}u)\right)\,dt=
\int\limits_0^\infty\left(e^{-\alpha t}\frac
d{dt}u(\gamma(t))-e^{-\alpha t}\alpha
u(\gamma(t))\right)\,dt=\int\limits_0^\infty \frac
d{dt}\left(e^{-\alpha t}u(\gamma(t))\right)\,dt=-u(\gamma(0)).$

Standard Hamilton--Jacobi equation corresponds to $\alpha=0$. As
we know, in general, this equation does not have smooth solutions,
it has only generalized ones. It would be very interesting to
study transformations of the solutions for the parameter $\alpha$
running from the indicated by Theorem~2 ``smooth area" to 0.

\medskip
\noindent{\bf Remark 2.}  The Hamiltonian $H$ is the energy of a
natural mechanical system on the Riemannian manifold and the
discount factor $\alpha$ plays the role of a negative friction
coefficient. Moreover, the change of variables $z\mapsto -z,\ z\in
T^*_qM,\ q\in M,$ transforms $\mathfrak h_\alpha$ in $-\mathfrak
h_{-\alpha}$ that allows to apply our analysis of the flow
$e^{t\mathfrak h_\alpha}$ to the dissipative mechanical system
(with a positive friction coefficient). As a byproduct we obtain
the description of the subset of $T^*M$ filled by the bounded
trajectories in the case when the friction coefficient is greater
than certain critical value.

%Indeed, $\alpha=0$ is not available in the setting of Theorem~2
%since $R^H_{0_q}=\nabla_q^2U>0$ if $q$ realizes a local minimum of
%$U$. Here  In the setting of Theorem~1, $\alpha=0$ is available
%only if $U$ is strongly concave; then the Lagrangian is strongly
%convex and the result follows from \cite[Th.\,4.1]{AgCh}.

\section{Partial hyperbolicity}

We start to prove Theorem~2. The Levi Civita connection $\nabla$ on
$T^*M$ defines a smooth ``horizontal" vector distribution
$D=\bigcup\limits_{z\in T^*M}D_z$ on $T^*M$, where $D_z$ is the
subspace of $T_z(T^*M)$ filled by the velocities of the parallel
translations of $z$ along the curves in $M$. We denote:
$\Delta_z=T_z(T^*_{\pi(z)}M)$ and call $\Delta=\bigcup\limits_{z\in
T^*M}\Delta_z$ the vertical distribution. Then
$T_z(T^*M)=\Delta_z\oplus D_z,\ \forall z\in T^*M$.

Note that both $\Delta_z$ and $D_z$ are Lagrange subspaces of the
symplectic space $T_z(T^*M)$. This is evident for $\Delta_z$; in
what concerns $D_z$, its property to be a Lagrange subspace is
just another way to say that the Levi Civita connection is
symmetric (i.\,e. torsion free). A vector distribution on a subset
of the symplectic manifold is called a {\it Lagrange vector
distribution} if its fibers are Lagrange subspaces of the tangent
spaces.

Let $w\in T_z(T^*M),\ w=w_{ver}+w_{hor}$, where $w_{ver}\in
\Delta_z,\ w_{hor}\in D_z$. We set $|w|=|w_{ver}|+|\pi_*w_{hor}|$
and thus define a canonical Riemannian structure on $T^*M$.

\begin{prop} Assume that $R^H_z<\frac{\alpha^2}4I,\ \forall z\in
B_H$. Then there exist continuous Lagrange distributions
$E^\pm=\bigcup\limits_{z\in\mathcal E_\alpha}E_z^\pm$ on $\mathcal
E_\alpha$ and positive constants $c_\pm,\,\varepsilon$ such that:
\begin{enumerate}
\item
$e^{t\mathfrak h_{\alpha}}_*E^\pm=E^\pm,\ \forall t\in\mathbb R$.
\item $|e^{t\mathfrak h_{\alpha}}_*w_-|\le
c_-e^{(\frac\alpha{2}-\varepsilon)t}|w_-|,\ |e^{t\mathfrak
h_{\alpha}}_*w_+|\ge c_+e^{(\frac\alpha{2}+\varepsilon)t}|w_+|,\
\forall\ w_\pm\in E^\pm,\ t\ge 0.$
\item $E^\pm_z\cap\Delta_z =0,\ \forall z\in\mathcal E_\alpha$.
\end{enumerate}
\end{prop}
{\bf Proof.} Recall that ${e^{t\mathfrak h_\alpha}}^*\sigma=e^{\alpha t}\sigma$;
hence $e^{t\mathfrak h_\alpha}_*$ transforms Lagrange subspaces
of the tangent spaces to $T^*M$ in the Lagrange subspaces. We define
the Jacobi curves and the curvature operators in the same way as for
the Hamiltonian flow (see Appendix). Namely, let $z\in T^*M$; then
$$
J_z^{\mathfrak h_\alpha}(t)=e^{-t\mathfrak h_\alpha}_*\Delta_{\zeta(t)},
\quad \zeta(t)=e^{t\mathfrak h_\alpha}(z),
$$
is a monotone decreasing curve in the Lagrange Grassmannian
$L\left(T_z(T^*M)\right)$ and $R_z^{\mathfrak
h_\alpha}=R_{J_z^{\mathfrak h_\alpha}}(0)$ is a self-adjoint
operator on $\Delta_z$. Obviously, $R_{J_z^{\mathfrak
h_\alpha}}(t)=e^{-t\mathfrak h_\alpha}_*R_{\zeta}^{\mathfrak
h_\alpha} e^{t\mathfrak h_\alpha}_*\bigr|_{J_z^{\mathfrak
h_\alpha}(t)}$ and $R_z^{\mathfrak h_0}=R^H_z$.

\begin{lemma} $R_z^{\mathfrak h_\alpha}=R_z^H-\frac{\alpha^2}4I$.
\end{lemma}
{\bf Proof.} We set $D^\alpha_z=\left(J_z^{\mathfrak
h_\alpha}\right)^\circ(0)\in L\left(T_z(T^*M)\right),\
D^\alpha=\bigcup\limits_{z\in T^*M}D^\alpha_z$, the canonical
Ehresmann connection associated to the flow $e^{t\mathfrak
h_\alpha}$. Then $D^0=D$, the Levi Civita connection.

Let $O_z$ be a neighborhood of $z$ in $T^*M$ and $v_\alpha(z'),\
z'\in O_z,$ a smooth ``vertical" vector field. According to the
terminology described in the Appendix, the field $v_\alpha$ is
parallel for the connection $D^\alpha$ along trajectories of the
flow $e^{t\mathfrak h_\alpha}$ if and only if $[\mathfrak
h_\alpha,v_\alpha]\in D^\alpha$. The connection $D^\alpha$ is
characterized by the following property: if $v_\alpha\in\Delta$
and $[\mathfrak h_\alpha,v_\alpha]\in D^\alpha$, then $[\mathfrak
h_\alpha,[\mathfrak h_\alpha,v_\alpha]]\in \Delta$. Moreover,
$$
[\mathfrak h_\alpha,[\mathfrak h_\alpha,v_\alpha]](z)=-R^{\mathfrak h_\alpha}_zv_\alpha(z).
$$

Let $v$ be a vertical vector field on $O_z$ that is parallel for the Levi Civita connection
$D$ along the trajectories of the flow $e^{t\mathfrak h_0}$ and is constant on the
vertical rays $\left(\mathbb R_+\zeta_t\right)\cap O_z $, where $\zeta_t=e^{t\mathfrak h_\alpha}(z)$ is any point of the contained 
in $O_z$ interval of the passing through $z$ trajectory of the flow $e^{t\mathfrak h_\alpha}$. Recall that
$\tau z'=\exp\left((\ln\tau)e\right)(z'),\ \tau>0,\ z'\in T^*M$, where $e$ is the Euler vector field. 
The linearity of the
Levi Civita connection implies: $\exp(se)_*D=D,\ \forall s\in\mathbb R$. Hence
$$
[e,v]=-v,\quad [e,[\mathfrak h_0,v]]=0.
$$

Now I claim that the vector field $v_\alpha$ defined by the formula
$v_\alpha(\zeta_t)=e^{\frac{\alpha}2t}v(\zeta_t)$  is parallel for the connection $D^\alpha$
along the curve $\zeta_t$. Indeed,
$$
[\mathfrak h_\alpha,v_\alpha](\zeta_t)=[\mathfrak h_0+\alpha e,v_\alpha](\zeta_t)=
e^{\frac{\alpha}2t}\left([\mathfrak h_\alpha,v](\zeta_t)-\frac{\alpha}2v(\zeta_t)]\right),
$$
$$
[\mathfrak h_\alpha,[\mathfrak
h_\alpha,v_\alpha]](\zeta_t)=e^{\frac{\alpha}2t}
\left(\frac{\alpha^2}4I-R_{\zeta_t}^H\right)v(\zeta_t)\in\Delta_{\zeta_t}.
$$
We have obtained the desired formula
$R^{\mathfrak h_\alpha}_z=R^H_z-\frac{\alpha^2}4I$ and the following characterization of
$D^\alpha$:
$$
D^\alpha_z=span\{[\mathfrak h_0,v](z)-\frac{\alpha}2v(z) :
v\in\mathfrak V_0\}, \eqno (4)
$$
where $\mathfrak V_0$ is the space of vertical vector fields on
$O_z$ that are parallel along the flow $e^{t\mathfrak h_0}$ for
the connection $D=D^0$ and are constant on the vertical rays
$\{\tau\zeta_t : \tau,t\in\mathbb R,\ \zeta_t,\tau\zeta_t\in O_z\}. \qquad
\square $

Next lemma is a simple generalization of the hyperbolicity test of
Lewowicz and Wojtkowski (see  \cite[Th.\,5.2]{Wo}). The proof is
almost verbal repetition of the proof from the cited paper and we
omit it.

Let $N$ be a Riemannian manifold, $X\in\mathrm{Vec}N$, and
$Q:TM\to\mathbb R$ a pseudo-Riemannian structure on $N$ (i.\,e. a
smooth field of nondegenerate quadratic forms $Q_z:T_zN\to\mathbb
R,\ z\in N$). Let $L_Xq:v\mapsto\frac
d{dt}Q(e^{tX}_*v)\bigr|_{t=0},\ v\in TM$, be the Lie derivative of
Q in the direction of X.

\begin{lemma} Let $\beta\in\mathbb R$ and $S\subset N$ be a compact
invariant subset of the flow $e^{tX},\ t\in\mathbb R$. If the
quadratic form $L_XQ-\beta Q$ is positive definite on $TN\bigr|_S$
then there exist invariant for the flow $e^{tX}$ continuous vector
distributions $E^\pm=\bigcup\limits_{z\in\mathcal
E_\alpha}E_z^\pm$ on $\mathcal E_\alpha$ and positive constants
$c_\pm,\varepsilon$ such that $(-1)^\pm Q\bigr|_{E^\pm}>0$ and
$$
|e^{tX}_*w_-|\le c_-e^{(\frac{\beta}2-\varepsilon)t}|w_-|,\
|e^{tX}_*w_+|\ge c_+e^{(\frac{\beta}2+\varepsilon)t}|w_+|,\
\forall\ w_\pm\in E^\pm,\ t\ge 0.
$$
\end{lemma}

Next lemma is a generalization of the earlier observation of Piotr
Przytycki (see \cite[Lemma\,2.1]{AgCh}). Let now $N$ be a
symplectic manifold endowed with the symplectic form $\sigma$,
$X\in\mathrm{Vec}M,\ \beta\in\mathbb R$, and
$L_X\sigma=\beta\sigma$. Let $\Lambda^i=\bigcup\limits_{z\in
N}\Lambda_z^i,\ \Lambda_z^i\in L(T_zN),\ i=0,1$, be two smooth
Lagrange distributions on $N$. We assume that
$\Lambda_z^0\cap\Lambda^1_z=0,\ \forall z\in N$.

Let $v\in T_zN$, then $v=v_0+v_1$, where $v_i\in\Lambda^i_z$. We
define the pseudo-Riemannian structure $Q_{\Lambda^0\Lambda^1}$ on
$N$ by the formula:
$$
Q_{\Lambda^0\Lambda^1}(v)=\sigma(v_0,v_1),\quad \forall v\in
T_zN,\ z\in N.
$$
We define the distributions $\Lambda^i(t)=e^{-tX}_*\Lambda^i$;
then $t\mapsto\Lambda^i_z(t)$ is a curve in the Lagrange Grassmannian
$L(T_zN),\ \forall z\in N,\ i=0,1.$

\begin{lemma} Let $S\subset N$ be a compact invariant subset of
the flow $e^{tX},\ t\in\mathbb R$. If the curve
$\Lambda^0_z(\cdot)$ is monotone decreasing and the curve
$\Lambda^1_z(\cdot)$ is monotone increasing, $\forall z\in S$,
then the form $Q_{\Lambda^0\Lambda^1}$ satisfies conditions of
Lemma~4.
\end{lemma}
{\bf Proof.} Let $V\in\mathrm{Vec}N,\ V\bigr|_S=V^0+V^1$, where
$V^i\in\Lambda^i$. We set $V(t)=e^{-tX}_*V,\ V^i(t)=e^{-tX}_*V^i,\
V(t)=V(t)_0+V(t)_1,\ V^i(t)=V^i(t)_0+V^i(t)_1$, where
$V^i(t)_j\in\Lambda^j,\ i,j=0,1$. We have:
$\left(L_XQ_{\Lambda^0\Lambda^1}\right)(V)= \frac
d{dt}\sigma\left(V(t)_1,V(t)_0\right)\bigr|_{t=0}$ and $\beta
Q_{\Lambda^0\Lambda^1}(V)=\beta\sigma(V^0,V^1)=\frac
d{dt}\sigma\left(V^1(t),V^0(t)\right)\bigr|_{t=0}.$ Then
$$
\sigma\left(V(t)_1,V(t)_0\right)=\sigma\left(V(t)_1,V(t)\right)=
\sigma\left(V(t)_1,V^1(t)\right)+\sigma\left(V(t)_1,V^0(t)\right)
$$
$$
=\sigma\left(V(t)_1,V^1(t)\right)-\sigma\left(V(t)_0,V^0(t)\right)
+\sigma\left(V(t),V^0(t)\right)
$$
$$
=\sigma\left(V(t)_1,V^1(t)\right)-\sigma\left(V(t)_0,V^0(t)\right)
+\sigma\left(V^1(t),V^0(t)\right).
$$
The differentiation with respect to $t$ at $t=0$ gives:
$$
\left(L_XQ_{\Lambda^0\Lambda^1}\right)(V)=
\underline{\dot\Lambda}^1(V^1)-\underline{\dot\Lambda}^0(V^0)+\beta
Q_{\Lambda^0\Lambda^1}(V)
$$
and the monotonicity assumptions imply that
$\underline{\dot\Lambda}^0(V^0)<0,\
\underline{\dot\Lambda}^1(V^1)>0. \qquad \square$

Note that the manifold $N=T^*M$, vector field $X=\mathfrak h_\alpha$, invariant
subset $S=\mathcal E_\alpha$, constant $\beta=\alpha$, and
distributions $\Lambda^0=\Delta,\ \Lambda^1=D^\alpha$ satisfy
conditions of Lemmas 4 and 5 if $R^H_z<\frac{\alpha^2}4I,\ \forall
z\in\mathcal E_\alpha$. Let $E^{\pm}$ be the invariant
distributions guaranteed by Lemma~4. To complete the proof of
Proposition~1 it remains to show that that $E^\pm_z$ are
transversal to $\Delta_z$ Lagrange subspaces.

We'll prove a little bit more. Namely, we are going to show that:
$$
E^\pm_z=\lim\limits_{t\to\mp\infty}J_z^{\mathfrak h_\alpha}(t)=
\lim\limits_{t\to\mp\infty}\left(J_z^{\mathfrak
h_\alpha}\right)^\circ(t),\quad \forall z\in\mathcal E_\alpha.
\eqno (5)
$$
Indeed, as we know (see Appendix) the limits
$\lim\limits_{t\to\pm\infty}J_z^{\mathfrak
h_\alpha}(t)=J_z^{\mathfrak h_\alpha}(\pm\infty)$ and
$\lim\limits_{t\to\pm\infty}\left(J_z^{\mathfrak
h_\alpha}\right)^\circ(t)=\left(J_z^{\mathfrak
h_\alpha}\right)^\circ(\pm\infty)$ do exist and are transversal to
$\Delta_z=J_z^{\mathfrak h_\alpha}(0)$. Moreover, $J_z^{\mathfrak
h_\alpha}(\pm\infty)$ and $\left(J_z^{\mathfrak
h_\alpha}\right)^\circ(\pm\infty)$ are invariant vector
distributions for the flow $e^{t\mathfrak h_\alpha}_*$, since
$J^{\mathfrak h_\alpha}(t+s)=e^{-t\mathfrak h_\alpha}J^{\mathfrak
h_\alpha}(s),\ \left(J^{\mathfrak
h_\alpha}\right)^\circ(t+s)=e^{-t\mathfrak
h_\alpha}\left(J^{\mathfrak h_\alpha}\right)^\circ(s)$.

Take a vector field $V=V^++V^-$, where $V^\pm\in E^\pm$. If
$V^\pm(z)\ne 0,\ \forall z\in\mathcal E_\alpha$, then the
component $e^{t\mathfrak h_\alpha}_*V^+\in E^+$ of the vector
$e^{t\mathfrak h_\alpha}_*V=e^{t\mathfrak
h_\alpha}_*V^++e^{t\mathfrak h_\alpha}_*V^-$ dominates as
$t\longrightarrow+\infty$ and the component $e^{t\mathfrak
h_\alpha}_*V^-$ dominates as $t\longrightarrow-\infty$ due to the
already proved estimates (see Lemma~4).

Therefore in order to prove (5) it is sufficient to show that
$J^{\mathfrak h_\alpha}(t)\cap E^\pm=\left(J^{\mathfrak
h_\alpha}\right)^\circ(t)\cap E^\pm=0$ for some (and hence for
all) $t\in\mathbb R$.

We have: $Q_{\Delta D^\alpha}\bigr|_{J^{\mathfrak
h_\alpha}(0)}=Q_{\Delta D^\alpha}\bigr|_\Delta=0$ and $\frac
d{dt}Q_{\Delta D^\alpha}\bigr|_{J^{\mathfrak h_\alpha}(0)}<0$.
Hence $Q_{\Delta D^\alpha}\bigr|_{J^{\mathfrak h_\alpha}(t)}<0$
for small positive $t$ and $Q_{\Delta
D^\alpha}\bigr|_{J^{\mathfrak h_\alpha}(t)}>0$ for small negative
$t$. On the other hand, $Q\bigr|_{E^+}>0$ and $Q\bigr|_{E^-}<0$.
It follows that $J^{\mathfrak h_\alpha}(t)\cap E^+=J^{\mathfrak
h_\alpha}(t)\cap E^-=0$. Similarly $\left(J^{\mathfrak
h_\alpha}\right)^\circ(t)\cap E^+=\left(J^{\mathfrak
h_\alpha}\right)^\circ(t)\cap E^-=0. \qquad \square$

\begin{corollary} Under conditions of Proposition~1, $\mathfrak
h_\alpha(z)\in E^-_z,\ \forall z\in\mathcal E_\alpha$.
\end{corollary}
{\bf Proof.} Let $\mathfrak h_\alpha=\mathfrak
h_\alpha^++\mathfrak h_\alpha^-$, where $\mathfrak
h_\alpha^\pm(z)\in E_z^\pm,\ \forall z\in\mathcal E_\alpha$. Then
the length of the vectors $e^{t\mathfrak h_\alpha}_*\mathfrak
h_\alpha^\pm(z)=\mathfrak h_\alpha^\pm(e^{t\mathfrak
h_\alpha}(z))$ is uniformly bounded as $t$ tends to $+\infty$. On
the other hand, $|e^{t\mathfrak h_\alpha}_*\mathfrak
h_\alpha^+(z)|\ge c_+e^{(\frac\alpha{2}+\varepsilon)t}|\mathfrak
h_\alpha^+(z)|$ for all $t\ge 0$. Hence $\mathfrak
h_\alpha^+(z)=0,\ \forall z\in\mathcal E_\alpha.\qquad \square$

\section{Optimal synthesis}

Now we are going to consider variational problems with finite
horizons and ``free endpoints" and then to study the limit as the
horizon tends to infinity. Namely, we study the functionals
$$
\mathfrak I^\tau_\alpha:\gamma\mapsto \int\limits_0^\tau e^{-\alpha
t}\left(\frac 12|\dot\gamma(t)|^2-U(\gamma(t))\right)\,dt
$$
on the spaces $\Omega_q^{\alpha,\tau}=\{\gamma\in
H^1([0,\tau];M):\gamma(0)=q\}.$ Compactness of $M$ implies the
existence of minimizers that are critical points of $\mathfrak
I_\alpha^\tau$ on $\Omega_q^{\alpha,\tau}$. These critical points
are solutions $\gamma$ of the Euler--Lagrange equations such that
the {\it transversality condition} $\dot\gamma(\tau)=0$ is
satisfied. In other words, critical points are projections to $M$ of
solutions $\zeta$ to equation (3) such that $\zeta(\tau)$ belongs to
the zero section of $T^*M$. Note that $\zeta(t)\in B_H,\ \forall
t\in[0,\tau]$; indeed, if a solution leaves $B_H$, then it never
comes back (see the prove of Lemma~2). In the sequel, we identify
$M$ with zero section of $T^*M$, i.\,e. $M=\{0_q:q\in M\}$, so that
$M\subset T^*M$.

\begin{prop} Assume that $R^H_z<\frac{\alpha^2}4I,\ \forall z\in
B_H$. Then $e^{-\tau\mathfrak h_\alpha}(M)$ is a smooth section of
the bundle $\pi:T^*M\to M,\ \forall \tau>0$; hence for any $q\in
M$ there exists a unique critical point of $\mathfrak
I_\alpha^\tau$ on $\Omega_q^{\alpha,\tau}$.
\end{prop}
{\bf Proof.} It is sufficient to show that the map $\pi\circ
e^{-t\mathfrak h_\alpha}\bigr|_M:M\to M$ has no critical points,
$\forall t>0$. Indeed, in this case the maps $\pi\circ
e^{-t\mathfrak h_\alpha}\bigr|_M$ must be coverings of $M$ and
$\pi\bigr|_M=id$; hence $\pi\circ e^{-t\mathfrak h_\alpha}\bigr|_M$
are actually diffeomorphisms.

So we have to show that $e^{-t\mathfrak
h_\alpha}_*T_{0_q}M\cap\Delta_{e^{t\mathfrak h_\alpha}(0_q)}=0$ for
any $q\in M$. To this end (and for further limiting procedure) we
introduce the Lagrange distribution
$\Delta^\alpha=\bigcup\limits_{z\in T^*M}\Delta^\alpha_z\subset
T(T^*M)$ as follows:
$$
\Delta^\alpha_z=span\left\{\left(1+\frac{\alpha}2\right)v-[\mathfrak
h_0,v] : v\in\mathfrak V_0\right\},\quad z\in T^*M,
$$
where $\mathfrak V_0$ has the same meaning as in (4). Consider the
splitting $T(T^*M)=\Delta^\alpha\oplus D^\alpha$. Symplectic form
$\sigma$ defines the nondegenerate pairing of the subspaces
$\Delta^\alpha_z$ and $D^\alpha_z$ that gives the identification
$\Delta^\alpha=(D^\alpha)^*$.

Then any transversal to $\Delta^\alpha_z$ Lagrange subspace
$\Lambda\subset T_z(T^*M)$ is identified with the graph of a
self-adjoint linear map from $D^\alpha_z$ to
$\Delta^\alpha_z=(D^\alpha)^*$ and thus with a quadratic form
$Q_\Lambda$ on $D^\alpha_z$.

\begin{lemma} We have:
$ Q_{D^\alpha_z}=0,\ Q_{\Delta_z}>0,\
Q_{D_z}=\frac{\alpha}{\alpha+2}Q_{\Delta_z},\ \forall z\in T^*M. $
In particular, $Q_{D^\alpha}<Q_D<Q_\Delta$.
\end{lemma}
{\bf Proof.} The equality $Q_{D^\alpha_z}=0$ is obvious. Let $\xi\in
D^\alpha_z$, then $\xi=[\mathfrak h_0,v](z)-\frac{\alpha}2v(z)$ for
some $v\in\mathfrak V_0$ (see (4)). Assume that $\xi\ne 0$, then
$v(z)\ne 0$. Moreover, $v(z)\in\Delta_z$ and $
v(z)=\left((1+\frac{\alpha}2)v-[\mathfrak
h_0,v]\right)+\left([\mathfrak h_0,v]-\frac{\alpha}2v\right), $
where $\left((1+\frac{\alpha}2)v-[\mathfrak
h_0,v]\right)\in\Delta^\alpha,\ \left([\mathfrak
h_0,v]-\frac{\alpha}2v\right)\in D^\alpha$. Then
$$
Q_{\Delta_z}(\xi)=\sigma_z\left(\left(1+\frac{\alpha}2\right)v-[\mathfrak
h_0,v]\,,\ [\mathfrak
h_0,v]-\frac{\alpha}2v\right)=\sigma_z([\mathfrak h_0,v],v)>0.
$$
Similarly, $Q_{D_z}(\xi)=\frac{\alpha}{\alpha+2}\sigma_z([\mathfrak
h_0,v],v).\qquad \square$

Note that $T_{0_q}M=D_{0_q},\ \forall q\in M$, since the Levi
Civita connection $D$ is a linear connection. Let
$z=e^{-\tau\mathfrak h_\alpha}(0_q),\ \zeta(t)=e^{t\mathfrak
h_\alpha}(z)$; then $t\mapsto J_z^{\mathfrak
h_\alpha}(t)=e^{-\tau\mathfrak h_\alpha}_*\Delta_{\zeta(t)}$ is a
monotone decreasing and $t\mapsto\left(J_z^{\mathfrak
h_\alpha}\right)^\circ(t)=e^{-\tau\mathfrak
h_\alpha}_*D^\alpha_{\zeta(t)}$ a monotone increasing curves in
the Lagrange Grassmannian $L\left(T_z(T^*M)\right)$. We'll use
simplified notations:
$$
\Delta(t)=e^{-t\mathfrak h_\alpha}_*\Delta_{\zeta(t)},\quad
D^\alpha(t)=e^{-t\mathfrak h_\alpha}_*D^\alpha_{\zeta(t)},\quad
D(t)=e^{-t\mathfrak h_\alpha}_*D_{\zeta(t)}.
$$
Then $t\mapsto Q_{\Delta(t)}$ is a strongly monotone decreasing and
$t\mapsto Q_{D^\alpha(t)}$ a strongly monotone increasing families
of quadratic forms. Moreover, $Q_{\Delta(t)}-Q_{D(t)}$ and
$Q_{D(t)}-Q_{D^\alpha(t)}$ are nondegenerate quadratic forms since
$D(t)$ is transversal to $\Delta(t)$ and $D^\alpha(t)$. It follows
that
$$
Q_{D^\alpha(0)}<Q_{D^\alpha(t)}<Q_{D(t)}<Q_{\Delta(t)}<Q_{\Delta(0)},
\quad\forall t>0. \eqno (6)
$$
In particular, $D(\tau)=e^{-\tau\mathfrak h_\alpha}_*D_{0_q}$ is
transversal to $\Delta(0)=\Delta_z$. Proposition~2 is proved.$\qquad
\square$

\medskip
Given $q\in M$ we denote by $\Phi_\tau(q)$ the value at $q$ of the
section $e^{-\tau\mathfrak h_\alpha}(M)$ of $T^*M$; in other words,
$\{\Phi_\tau(q)\}=e^{-\tau\mathfrak h_\alpha}(M)\cap T^*_qM$. Recall
that $\Phi_\tau(q)\in B_H,\ \forall q\in M,\ \tau\ge 0$. In
particular, $\tau\mapsto\Phi_\tau(q),\ \tau\ge 0,$ is a uniformly
bounded curve in $T^*M$.

\begin{lemma} Assume that
$z=\lim\limits_{k\to\infty}\Phi_{\tau_k}(q)$ for some
$\tau_k\longrightarrow+\infty$ as $k\longrightarrow\infty$. Then
$z\in\mathcal E_\alpha$ and
$E^-_z=\lim\limits_{k\to\infty}{\Phi_{\tau_k}}_*(T_qM)$.
\end{lemma}
{\bf Proof.} Let $\gamma(t)=\pi\circ e^{t\mathfrak h_\alpha}(z),\
\gamma_k(t)=\pi\circ e^{t\mathfrak
h_\alpha}\left(\Phi_{\tau_k}(q)\right)$. Then $\mathfrak
I^\tau_\alpha(\gamma_k)\le\mathfrak
I^\tau_\alpha(\gamma_k)(const)=\frac 1\alpha(e^{-\alpha
t}-1)U(q)<\frac 1\alpha |U(q)|$. Hence $\mathfrak
I^\tau_\alpha(\gamma)=\lim\limits_{k\to\infty}\mathfrak
I^\tau_\alpha(\gamma_k)\le\frac 1\alpha |U(q)|,\ \forall\tau\ge 0$.
We obtain that $\mathfrak I_\alpha(\gamma)\le \frac 1\alpha |U(q)|$.
In particular, $\gamma\in\Omega^\alpha_q$ and, according to Lemma~1,
$z\in\mathcal E_\alpha$.

The subspace ${\Phi_{\tau_k}}_*(T_qM)\subset
T_{\Phi_{\tau_k}(q)}(T^*M)$ is the tangent space to the submanifold
$e^{-\tau_k\mathfrak h_\alpha}(M)$, i.\,e.
${\Phi_{\tau_k}}_*(T_qM)=e^{-\tau_k\mathfrak
h_\alpha}_*\bigl(D_{0_{\gamma_k(\tau_k)}}\bigr).$ Now the statement
of the lemma follows from (5) and (6). $\qquad \square$

Recall that $\Phi_\tau(M)$ is an exact Lagrange submanifold of $T^*M$,
hence $s\bigr|_{\Phi_\tau(M)}$ is an exact 1-form, where $s$ is the
Liouville form on $T^*M$.
%Let $\varpi:\tilde M\to M$ be the
%universal covering of $M$, $\tilde\mathfrak
%h_\alpha\in\mathrm{Vec}(T^*\tilde M)$ the pullback of the vector
%field $\mathfrak h_\alpha$ and $\tilde\Phi_\tau:T^*\tilde M\to\tilde
%M$ the pullback of the section $\Phi_\tau$. Let $\tilde s$ be the
%Liouville form on $T^*\tilde M$, then $\tilde
%s\bigr|_{\tilde\Phi_\tau(\tilde M)}$ is an exact form.
In other words, $\Phi_\tau=da_\tau$ for a smooth scalar function
$a_\tau$ on $ M$. Lemma~7 together with statement {\it 3} of
Proposition~1 imply that second derivatives of $a_\tau$ are
uniformly bounded for all $\tau\ge 0$.

%Fix an open subset $\mathcal O\subset\tilde M$ such that the
%closer of $\mathcal O$ is compact and $\varpi(\mathcal O)=M$.
The functions $a_\tau$ are defined up to a constant and we may of
course assume that they are uniformly bounded on $M$.
Then there exists a sequence $\tau_k\longrightarrow+\infty$ as
$k\longrightarrow\infty$ and function $a_\infty\in
C^{1,\infty}(M)$ such that $a_{\tau_k}\longrightarrow
a_\infty$ as $k\longrightarrow\infty$ in $C^1$-topology.
%Moreover, $da_\infty$ is the pullback of a 1-form $\psi$ on $M$, where
Set $\psi(q)=d_qa_\infty$, then
$\psi(q)=\lim\limits_{k\to\infty}\Phi_{\tau_k}(q),\ \forall q\in
M$. We obtain that $\psi(q)\in\mathcal E_\alpha$ and
${\Phi_{\tau_k}}_*(T_qM)\longrightarrow E^-_{\psi(q)}$ as
$k\longrightarrow\infty,\ \forall q\in M$. Hence the function
$a_\infty$ is actually of class $C^2$; the submanifold
$\Psi\stackrel{def}{=}\{\psi(q):q\in M\}$ is contained in
$\mathcal E_\alpha$ and $T_z\Psi=E^-_z,\ \forall z\in\Psi$.

According to Corollary~3, $\mathfrak h_\alpha(z)\in E^-_z$, hence
$\Psi$ is an invariant exact Lagrange submanifold of the flow
$e^{t\mathfrak h_\alpha},\ t\in\mathbb R$. Moreover, $\Psi\supset
C_H$ since $C_H=e^{\tau\mathfrak h_\alpha}(C_H)\subset
e^{\tau\mathfrak h_\alpha}(M),\ \forall\tau\ge 0$.

If $U$ is a Morse function, then any fixed point $z=0_q\in C_H$ of
the flow $e^{t\mathfrak h_\alpha}$ is hyperbolic with real
eigenvalues. This is immediately seen after the diagonalization of
the Hessian of $U$ in the critical point $q$ by the orthogonal
transformation of $T_qM$. The stable subspace of the linearization
of $\mathfrak h_\alpha$ at $z$ is contained in $E^-_z=T_z\Psi$.
Hence the stable submanifold of the hyperbolic equilibrium $z$ is
contained in $\Psi\subset\mathcal E_\alpha$. On the other hand,
$\mathcal E_\alpha$ is the union of the stable submanifolds of all
equilibria $z\in C_H$ (see Lemmas 1,\,2). We obtain that
$\Psi=\mathcal E_\alpha$.

Figure 5 illustrates the structure of $\mathcal E_\alpha$ near an
equilibrium point $z$. The stable subspace of the linearized system
is always contained in $E^-_z=T_z\mathcal E_\alpha$, while the
unstable subspace splits in the ``less unstable" part that is
contained in $E^-_z$ and the ``more unstable" part that is equal to
$E^+_z$.

\begin{figure}
\hspace{3cm} \begin{minipage}[b]{6.5cm}
   \centering
   \includegraphics[width=7cm]{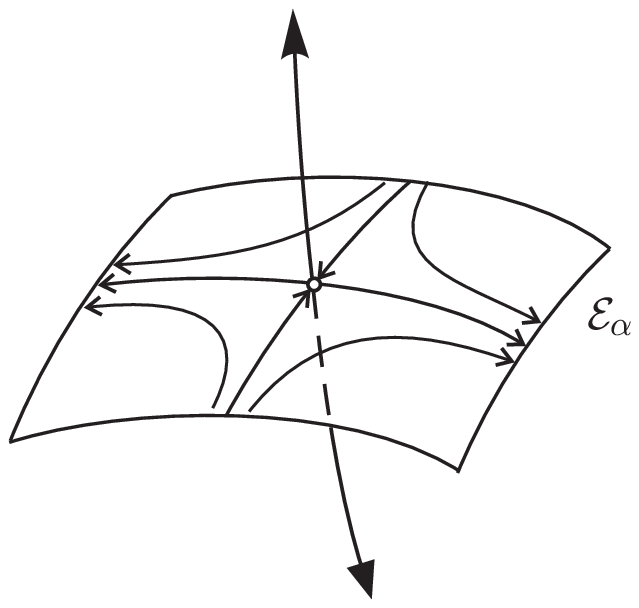}
   \caption{}
 \end{minipage}
 \end{figure}

We have proved all good properties of $\mathcal E_\alpha$ stated
in Theorem~2. It remains only to check that the curves $t\mapsto
\pi\circ e^{t\mathfrak h_\alpha}(z),\ z\in\mathcal E_\alpha$, are
minimizers. To do that, we use the classical ``fields of
extremals" method.
%Let $\tilde\mathcal E_\alpha\subset T^*\tilde
%M$ be the pullback of $\mathcal E_\alpha\subset T^*M$.
We set:
$$
\mathcal L=\left\{(e^{-\alpha t}z,t): z\in\mathcal E_\alpha,\
t\ge 0\right\}\subset T^*M\times\mathbb R_+.
$$
Then $\mathcal L$ is an $n$-dimensional submanifold of $T^*
M\times\mathbb R_+$ and 1-form \linebreak $\bigl(
s-e^{-\alpha t}Hdt\bigr)\Bigr|_{\mathcal L}$ is exact.
% here $\tilde H:T^*\tilde M\to\mathbb R$ is the pullback of the
%Hamiltonian $H$.

Given $q\in M$, there exists a unique $z\in\mathcal E_\alpha$
such that $\pi(z)=q$. We set $\gamma(t)=\pi\circ e^{t\mathfrak
h_\alpha}(z)$ and we have to prove that $\mathfrak
I_\alpha(\gamma)<\mathfrak I_\alpha(\varrho),\
\forall\varrho\in\Omega_q^\alpha$ such that $\varrho\ne\gamma$.

%Take $\tilde q\in\tilde M$ such that $\varpi(\tilde q)=q$ and the
%curves $\tilde\gamma,\ \tilde\varrho$ in $\tilde M$ such that
%$\tilde\gamma(0)=\tilde\varrho(0)=\tilde q,\
%\varpi\circ\tilde\gamma=\gamma,\ \varpi\circ\tilde\varrho=\varrho$;
%these curves do exist and are uniquely defined since \linebreak
%$\varpi:\tilde M\to M$ is a covering. Obviously, $\mathfrak
%I_\alpha(\gamma)= \tilde\mathfrak{I}_\alpha(\tilde\gamma)$ and
%$\mathfrak
%I_\alpha(\varrho)=\tilde\mathfrak{I}_\alpha(\tilde\varrho)$, where
%$\tilde\mathfrak{I}_\alpha$ is the pullback of $\mathfrak I_\alpha$.

Let $$ \zeta(t)\in\mathcal E_\alpha\cap
T_{\gamma(t)} M,\qquad \eta(t)\in\mathcal
E_\alpha\cap T_{\varrho(t)} M
$$
be the lifts of $\gamma$ and $\varrho$ to
$\mathcal E_\alpha$ and
$$
\hat\zeta(t)=\left(e^{-\alpha t}\zeta(t),t\right),\qquad
\hat\eta(t)=\left(e^{-\alpha t}\eta(t),t\right)
$$
be the lifts of $\gamma(t)$ and $\varrho(t)$ to
$\mathcal L$. Then
$$
\int\limits_{\hat\zeta}( s-e^{-\alpha t} Hdt)=
\int\limits_0^\infty e^{-\alpha
t}\bigl(\langle\zeta(t),\dot\gamma(t)\rangle-H(\zeta(t))\bigr)\,dt=\mathfrak
I_\alpha(\gamma)
$$
and
$$
\int\limits_{\hat\eta}( s-e^{-\alpha t} Hdt)=
\int\limits_0^\infty e^{-\alpha
t}\bigl(\langle\eta(t),\dot\varrho(t)\rangle-H(\eta(t))\bigr)\,dt<
\mathfrak I_\alpha(\varrho).
$$
Now, $\forall\tau\ge 0$, we define the curves
$\hat\zeta_\tau:[0,2\tau]\to\mathcal L$ and
$\hat\eta_\tau:[0,2\tau]\to\mathcal L$ by the formulas:
$$
\hat\zeta_\tau(t)=\left\{\begin{array}{cl} \hat\zeta(t),& 0\le t\le\tau,\\
\left(e^{-\alpha\tau}\zeta(2\tau-t),\tau\right),& \tau\le t\le 2\tau;\\
\end{array}\right.
$$
$$
\hat\eta_\tau(t)=\left\{\begin{array}{cl} \hat\eta(t),& 0\le t\le\tau,\\
\left(e^{-\alpha\tau}\eta(2\tau-t),\tau\right),& \tau\le t\le 2\tau.\\
\end{array}\right.
$$
The curves $\hat\zeta_\tau$ and $\hat\eta_\tau$ have common
endpoints, hence
$$
\int\limits_{\hat\zeta_\tau}(s-e^{-\alpha t}Hdt)=
\int\limits_{\hat\eta_\tau}(s-e^{-\alpha t}Hdt).
$$
We have:
$$
\int\limits_{\hat\eta_\tau}(s-e^{-\alpha t}Hdt)=
\int\limits_{\hat\eta|_{[0,\tau]}}(s-e^{-\alpha t}
Hdt)-
e^{-\alpha\tau}\int\limits_0^\tau\langle\eta(t),\dot\varrho(t)\rangle\,dt
$$
and
$$
\Bigl|e^{-\alpha\tau}\int\limits_0^\tau\langle\eta(t),\dot\varrho(t)\rangle\,dt\Bigr|
<ce^{-\frac{\alpha}2\tau}\int\limits_0^\tau|e^{-\frac{\alpha}2t}\dot\varrho(t)|\,dt
$$
$$
\le ce^{-\frac{\alpha}2\tau}\sqrt{\tau}\left(\int\limits_0^\infty
e^{-\alpha t}|\dot\varrho(t)|^2\,dt\right)^{\frac 12},
$$
where $c=\max\{|z|:z\in\mathcal E_\alpha\}$. Hence
$$
\lim\limits_{\tau\to+\infty}\int\limits_{\hat\eta_\tau}(
s-e^{-\alpha t}H dt)=\int\limits_{\hat\eta}(
s-e^{-\alpha t}H dt)
$$ and similarly
$$
\lim\limits_{\tau\to+\infty}\int\limits_{\hat\zeta_\tau}(
s-e^{-\alpha t}H dt)=\int\limits_{\hat\zeta}(
s-e^{-\alpha t}H dt).
$$
Summing up, we obtain:
$$
\mathfrak I_\alpha(\gamma)=\int\limits_{\hat\zeta}(
s-e^{-\alpha t}H dt)=\int\limits_{\hat\eta}(
s-e^{-\alpha t}H dt)< \mathfrak I_\alpha(\varrho). \eqno
\square
$$

%We continue the proof of Theorem~2.

%\begin{lemma} Let $q_0\in M$ be a point where $U$ attains
%maximum, i.\,e. $U(q_0)=\max U$, and let $\mathcal E^0_\alpha$ be the
%connected component of $\mathcal E_\alpha$, which contains
%$0_{q_0}$ (the origin of $T^*_{q_0}$). Then $\mathcal E^0_\alpha$ is a $C^1$-submanifold of
%$T^*M$ and $T_z\mathcal E^0_\alpha=E^-_z,\ \forall z\in \mathcal
%E^0_\alpha$.
%\end{lemma}
%{\bf Sketch of proof.}

%\begin{corollary} $\pi\bigr|_{\mathcal E^0_\alpha}$ is
%a diffeomorphism of $\mathcal E^0_\alpha$ on $M$.
%\end{corollary}
%{\bf Proof.} Let $z\in\mathcal E^0_\alpha$, we have:
%$\ker\pi_*\bigr|_{T_z(T^*M)}=\Delta_z, T_z\mathcal
%E^0_\alpha=E^-_z$ and $\Delta_z\cap E^-_z=0$. Hence
%$\pi\bigr|_{\mathcal E^0_\alpha}$ is a local diffeomorphism.
%Moreover, $\mathcal E^0_\alpha$ is compact, hence
%$\pi\bigr|_{\mathcal E^0_\alpha}$ is a covering. On the other
%hand, $\mathcal E_\alpha\subset B_H$ and $B_H\cap
%T^*_{q_0}M=0_{q_0}$. It follows that
%our covering has only one leaf and is actually a
%diffeomorphism.\qquad $\square$

\section*{Appendix}

Here we collect some definitions and geometric facts that are used
in the paper; see \cite{Ag} and references therein for the
consistent presentation.

\subsection*{Monotone curves in the Lagrange Grassmannian}

Let $\Sigma$ be a $2n$-dimensional symplectic space endowed with a
symplectic form $\sigma$. A Lagrange subspace is an $n$-dimensional
subspace $\Lambda\subset\Sigma$ such that
$\sigma\bigr|_{\Lambda}=0$. The set of all Lagrange subspaces forms
a compact $\frac{n(n+1)}2$-dimensional manifold that is called
Lagrange Grassmannian and is denoted $L(\Sigma)$. Symplectic group
acts transitively on the Lagrange Grassmannian; moreover, symplectic
group acts transitively on the pairs of transversal Lagrange
subspaces.

Let $\Pi\in L(\Sigma)$, the symplectic form defines a
non-degenerate pairing of $\Pi$ and $\Sigma/\Pi$: the scalar
product of $x\in\Pi$ and the residue class $y+\Pi$ is equal to
$\sigma(x,y)$. We obtain a natural identification
$\Sigma/\Pi=\Pi^*$. Now set:
$$
\Pi^\pitchfork=\{\Lambda\in L(\Sigma): \Pi\cap\Lambda=0\};
$$
then $\Pi^\pitchfork$ is an affine space over the vector space
$\mathrm{Sym}(\Sigma/\Pi)$ of linear self-adjoint maps from
$\Sigma/\Pi$ to $\Pi=(\Sigma/\Pi)^*$ (or, in other words, over the
space of quadratic forms on $\Sigma/\Pi$). Indeed, the sum of
$\Lambda\in\Pi^\pitchfork$ and $S\in\mathrm{Sym}(\Sigma/\Pi)$ is
defined as follows:
$$
\Lambda+S=\{S(y+\Pi)+y : y\in\Lambda\}\in L(\Sigma).
$$

An affine space is a vector space ``with no origin". Let us take
$\Delta\in\Pi^\pitchfork$ and order $\Delta$ to be the origin.
Then $\Pi^\pitchfork$ turns into $\mathrm{Sym}(\Sigma/\Pi)$.
Moreover, obvious isomorphism $\Sigma/\Pi\cong\Delta$ induces the
isomorphism of $\mathrm{Sym}(\Sigma/\Pi)$ with
$\mathrm{Sym}\Delta$ and the isomorphism of $\Pi$ with $\Delta^*$.
This makes $\Pi^\pitchfork$ a coordinate chart of the manifold
$L(\Sigma)$ coordinatized by $\mathrm{Sym}(\Delta)$.

Given $\Lambda_i\in\Pi^\pitchfork,\ i=0,1$, let
$Q_{\Lambda_i}\in\mathrm{Sym}(\Delta)$ be the coordinate
presentation of $\Lambda_i$; then
$\dim(\Lambda_0\cap\Lambda_1)=\dim\ker(Q_{\Lambda_0}-Q_{\Lambda_1})$.

Let $V\in T_{\Lambda_0}L(\Sigma)$. To the tangent vector $V$ we
associate a quadratic form $\underline{V}$ on $\Lambda_0$ (or, in
other words, a self-adjoint map from $\Lambda_0$ to
$\Lambda_0^*=\Sigma/\Lambda_0$) as follows: take a smooth curve
$\Lambda_t\in L(\Sigma)$ such that $\dot\Lambda_0=V$ and a smooth
curve $x_t\in\Lambda_t$ and set $\underline{V}(x_0)=\sigma(x_0,\dot
x_0)$. Then $V\mapsto\underline{V}$ is a natural isomorphism of
$T_{\Lambda_0}L(\Sigma)$ and $\mathrm{Sym}(\Lambda_0)$.

A curve $t\mapsto\Lambda_t$ is {\it regular} if
$\underline{\dot\Lambda}_t$ is a non-degenerate quadratic form. A
regular curve is {\it monotone increasing} ({\it monotone
decreasing}) if $\underline{\dot\Lambda}_t$ is positive definite
(negative definite).

Let $\Pi\cap\Lambda_t=0$, then  $\Lambda_t$ belongs to the affine
space $\Pi^\pitchfork$ and the derivative $\dot\Lambda_t$ belongs to
the vector space $\mathrm{Sym}(\Sigma/\Pi)$. Obvious isomorphism
$\Sigma/\Pi\cong\Lambda_t$ induces the isomorphism
$\mathrm{Sym}(\Sigma/\Pi)\cong\mathrm{Sym}\Lambda_t$, which
transforms $\dot\Lambda_t$ in $\underline{\dot\Lambda}_t$. In
particular, a monotone increasing (decreasing) curve is presented by
a strongly monotone increasing (decreasing) family of quadratic
forms in any coordinate chart $\Pi^\pitchfork$.

Let $t\mapsto\Lambda_t$ be a regular curve in $L(\Sigma)$ and
$\tau\in\mathbb R$; then $\Lambda_t\in\Lambda_\tau^\pitchfork$ for
all $t$ that are sufficiently close and not equal to $\tau$. We can
treat $t\mapsto\Lambda_t$ as a curve in the affine space
$\Lambda_\tau^\pitchfork$ with a singularity at $t=\tau$. This
singularity is actually a simple pole. We can write the Laurent
expansion of $\Lambda_t$ at $t=\tau$ in this affine space (that is
an affine space over the vector space
$\mathrm{Sym}(\Sigma/\Lambda_\tau)$). All terms of the Laurent
expansion but the free term are elements of the vector space
$\mathrm{Sym}(\Sigma/\Lambda_\tau)$, while the free term is an
element of the affine space $\Lambda_\tau^\pitchfork$ itself. We
denote this free term by $\Lambda_\tau^\circ$ and call the curve
$\tau\mapsto\Lambda_\tau^\circ$ the {\it derivative curve} of the
curve $\Lambda_\cdot$.

Given $t\in\mathbb R$, the derivative curve defines a splitting:
$\Sigma=\Lambda_t\oplus\Lambda^\circ_t$, which induces the
identifications
$$
\Lambda^\circ_t=\Sigma/\Lambda_t=\Lambda_t^*,\qquad
\Lambda_t=\Sigma/\Lambda^\circ_t={\Lambda^\circ_t}^*.
$$
Then
$$
\underline{\dot\Lambda}_t:\Lambda_t\to\Lambda_t^\circ,\qquad
\underline{\dot\Lambda}^\circ_t:\Lambda^\circ_t\to\Lambda_t
$$
are self-adjoint linear maps. The operator
$R_\Lambda(t):\Lambda_t\to\Lambda_t$ defined by the formula
$$R_\Lambda(t)=-\underline{\dot\Lambda}_t^\circ\circ\underline{\dot\Lambda}_t$$
is the {\it curvature operator} of the curve $\Lambda_\cdot$.

Assume that $t\mapsto\Lambda_t$ is a monotone curve, then
$|\underline{\dot\Lambda}_t|$ is a positive definite quadratic
form on $\Lambda_t$ and the curvature operator $R_\Lambda(t)$ is a
symmetric operator with respect to the Euclidean structure defined
by this quadratic form. In particular, the operator $R_\Lambda(t)$
is diagonalizable and all its eigenvalues are real. We say that
the monotone curve $\Lambda_\cdot$ {\it has a positive (negative)
curvature} if all eigenvalues of $R_\Lambda(t)$ are positive
(negative) and uniformly separated from 0. If a monotone curve
$\Lambda_\cdot$ has a positive (negative) curvature, then the
curve $\Lambda^\circ_\cdot$ is also monotone and the direction of
monotonicity of $\Lambda^\circ_\cdot$ coincides with (is opposite
to) the direction of monotonicity of $\Lambda_\cdot$.

Let $t\mapsto\Lambda_t$ be a regular curve and $t\mapsto\Delta_t$
another curve in $L(\Sigma)$ such that $\Lambda_t\cap\Delta_t=0,\
\forall t\in\mathbb R$. We may treat $\{(t,\Lambda_t):t\in\mathbb
R\}\subset\mathbb R\times\Sigma$ and $\{(t,\Delta_t):t\in\mathbb
R\}\subset\mathbb R\times\Sigma$ as subbundles of the trivial vector
bundle $\mathbb R\times\Sigma$; these subbundles define a splitting
of the trivial bundle. We say that the ``section" $x_t\in\Lambda_t,\
t\in\mathbb R,$ is {\it parallel for the splitting}
$\Sigma=\Lambda_t\oplus\Delta_t,\ t\in\mathbb R,$ if $\dot
x_t\in\Delta_t,\ \forall t\in\mathbb R$.

Canonical splitting $\Sigma=\Lambda_t\oplus\Lambda_t^\circ$ can be
characterized as follows: the relations $x_t\in\Lambda_t,\ \dot
x_t\in\Delta_t,\ \forall t\in\mathbb R,$ imply the relation $\ddot
x_t\in\Lambda_t$ if and only if $\Delta_t=\Lambda^\circ_t,\ \forall
t\in\mathbb R$; moreover, $\ddot x_t=-R_\Lambda(t)x_t$ in the case
of the canonical splitting.

\medskip\noindent{\bf Theorem}
{\it Let $\Lambda_t,\ t\in\mathbb R$, be a monotone curve in
$L(\Sigma)$. If The curve $\Lambda_\cdot$ has negative curvature,
then there exist} $
\lim\limits_{t\to\pm\infty}\Lambda_t=\Lambda_{\pm\infty}=
\lim\limits_{t\to\pm\infty}\Lambda^\circ_t$ {\it and} \linebreak
$\Lambda_t\cap\Lambda_\tau=0,\ \forall\, -\infty\le
t<\tau\le+\infty.$

\medskip
The existence of the limits and the fact that $\Lambda_t,\
t\in\mathbb R,$ are mutually transversal can be easily explained.
Recall that symplectic group acts transitively on the set of pairs
of transversal Lagrange subspaces and that coordinate presentation
of a Lagrange subspace $\Lambda\in L(\Sigma)$ is a quadratic form
$Q_\Lambda$. Assume that the curve $t\mapsto\Lambda_t$ is monotone
decreasing and the curve $t\mapsto \Lambda^\circ_t$ is monotone
increasing (the opposite ``increasing--decreasing" case is treated
similarly). We can always find an appropriate coordinate chart
such that $Q_{\Lambda_0}-Q_{\Lambda^\circ_0}>0$. Moreover,
$t\mapsto Q_{\Lambda_t}$ is a strictly monotone decreasing family
of quadratic forms, while $t\mapsto Q_{\Lambda^0_t}$ is a strictly
monotone increasing family, and
$Q_{\Lambda_t}-Q_{\Lambda^\circ_t}$ are non-degenerate forms.
Hence $\Lambda_t$ and $\Lambda^\circ_t$ never leave our chart for
positive $t$ and there exist
$\lim\limits_{t\to+\infty}Q_{\Lambda_t}\ge\lim\limits_{t\to+\infty}Q_{\Lambda^\circ_t}$.
A more carefull analysis shows that these two limits coincide and
convergence to the common limit has an exponential rate. The
limiting procedure as $t\longrightarrow-\infty$ is performed
similarly: we simply take a coordinate chart such that
$Q_{\Lambda_0}-Q_{\Lambda^\circ_0}<0$.

\subsection*{Jacobi curves}
Let $M$ be a smooth $n$-dimensional manifold and $T^*M$ its
cotangent bundle equipped with the standard symplectic structure.
Given $q\in M$ and $z\in T_q^*M$ we set $\Sigma_z=T_z(T^*M),\
\Delta_z=T_z(T^*_qM)$; then $\Sigma_z$ is a symplectic space and
$\Delta_z$ is a Lagrange subspace of this symplectic space.

Let $h:T^*M\to\mathbb R$ be a smooth (Hamiltonian) function, $\vec
h\in\mathrm{Vec}(T^*M)$ the associated to $h$ Hamiltonian vector
field, and $t\mapsto e^{t\vec h}$ the generated by $\vec h$
Hamiltonian flow on $T^*M$. {\sl Jacobi curve} $J_z^{\vec h}(t)$ is
a curve in the Lagrange Grassmannian $L(\Sigma_z)$ defined by the
formula:
$$
J_z^{\vec h}(t)=e^{-t\vec h}_*\Delta_{e^{t\vec h}(z)}, \quad
t\in\mathbb R.
$$

We list some easily derived from the definition basic properties
of the Jacobi curves. Let $\zeta(t)=e^{t\vec h}(z),\ \zeta(t)\in
T^*_{\gamma(t)}M$; then quadratic form $\underline{\dot J_z^{\vec
h}}(0)$ on $J_z^{\vec h}(0)=\Delta_z$ is equal to
$-d^2_z\left(h\bigr|_{T^*_qM}\right)$ and quadratic form
$\underline{\dot J_z^{\vec h}}(t)$ on $J_z^{\vec h}(t)$ is
obtained from
$-d^2_{\zeta(t)}\left(h\bigr|_{T^*_{\gamma(t)}M}\right)$ by the
linear change of variables \linebreak $e^{-t\vec
h}_*:\Delta_{\zeta(t)}\to J_z^{\vec h}(t)$. It follows that Jacobi
curve $J_z^{\vec h}(t)$ is monotone decreasing (increasing),
$\forall z\in T^*M$, if and only if $h\bigr|_{T^*_qM}$ is strongly
convex (concave), $\forall q\in M$. Other properties:
$$
{J_z^{\vec h}}^\circ(t)= e^{-t\vec h}_*{\left(J_{\zeta(t)}^{\vec
h}\right)}^\circ(0), \quad R_{J_z^{\vec h}}(t)=e^{-t\vec
h}_*R_{J_{\zeta(t)}^{\vec h}}(0)e^{t\vec h}_*\bigr|_{J_z^{\vec
h}(t)}.
$$

Lagrange distribution ${J_z^{\vec h}}^\circ(0),\ z\in T^*M$, on
$T^*M$ is called the {\it canonical connection} associated to $\vec
h$ and linear operator $R_{J_z^{\vec h}}(0):\Delta_z\to\Delta_z$ is
called the {\it curvature operator} of $\vec h$ at $z\in T^*M.$ We
will use simplified notations:
$$
\Delta^h_z={J_z^{\vec h}}^\circ(0), \qquad R^h_z=R_{J_z^{\vec
h}}(0);
$$
then $\Sigma_z=\Delta_z\oplus\Delta^h_z,\ z\in T^*M$, is the
associated to $\vec h$ {\it canonical splitting} of the vector
bundle $T(T^*M)$.

Any vector field $f\in\mathrm{Vec}(T^*M)$ splits in the {\it
vertical} and {\it horizontal} parts as follows:
$f=f_{ver}+f_{hor}$, where $f_{ver}(z)\in\Delta_z,\
f_{hor}(z)\in\Delta^h_z,\ \forall z\in T^*M$. Now let
$v\in\mathrm{Vec}(T^*M)$ be a vertical vector field, i.\,e.
$v_{hor}=0$; we say that the field $v$ is {\it parallel for the
connection} $\Delta^h$ along trajectories of the flow $e^{tf}$ if
$[f,v]_{ver}=0$.

Horizontal vector fields and parallel vertical vector fields can be
defined for any Ehresmann connection (i.\,e. for any vector
distribution $\mathfrak D$ on $T^*M$ such that
$\Sigma_z=\Delta_z\oplus\mathfrak D_z$). The associated to $\vec h$
canonical connection is characterized by the following property: if
$v_{hor}=0$ and $[\vec h,v]_{ver}=0$, then $[\vec h,[\vec
h,v]]_{hor}=0$.

Finally, for any vertical vector field $v$ and any $z\in T^*M$ we
have:
$$
R^h_zv(z)=-[\vec h,[\vec h,v]_{ver}]_{hor}(z).
$$

\medskip\noindent{\bf Example.} Let $M$ be a Riemannian manifold and
$H:T^*M\to\mathbb R$ the energy function of a natural mechanical
system on $M$, i.\,e. $H(z)=\frac 12|z|^2+U(q)$, \linebreak
$\forall q\in M,\ z\in T^*M$, where $U:M\to\mathbb R$ is the
potential energy. Then $\Delta^H$ is actually standard Levi Civita
connection $\nabla$ rewritten as an Ehresmann connection on
$T^*M$. More precisely, $\Delta^h_z$ is the subspace of $\Sigma_z$
filled by the velocities of the parallel translations of the
covector $z$ along curves in $M$. Moreover, Riemannian structure
gives the identification $T^*M\cong TM$. Combining this
identification with the identification $T^*_qM\cong
T_z(T^*_qM)=\Delta_z$ of the vector space $T_q^*M$ with its
tangent space we obtain the explicit formula for the curvature
operator $R^H_z$ of the Hamiltonian field $\vec H$ at $z\in T^*M$:
$$
R^H_z=\mathfrak R(\cdot,z)z+\nabla^2_qU,
$$
where $\mathfrak R$ is the Riemannian curvature and $\nabla^2$ is
the second covariant derivative.

\end{document}